\documentclass[12pt,oneside,english]{amsart}
\usepackage[T1]{fontenc}
\usepackage[latin1]{inputenc}
\usepackage{amssymb}

\makeatletter
 \theoremstyle{plain}    
 \newtheorem{thm}{Theorem}[section]
 \numberwithin{equation}{section} 
 \numberwithin{figure}{section} 
 \theoremstyle{plain}
 \theoremstyle{plain}    
 \newtheorem{lem}[thm]{Lemma} 
 \theoremstyle{definition}
  \newtheorem{condition}[thm]{Condition}
 \theoremstyle{definition}
 \newtheorem{defn}[thm]{Definition}
 \theoremstyle{plain}    
 \newtheorem{prop}[thm]{Proposition} 
 \theoremstyle{plain}    
 \newtheorem{cor}[thm]{Corollary} 
 \theoremstyle{remark}    
 \newtheorem*{note*}{Note} 

\usepackage{graphicx}

\usepackage{babel}
\makeatother
\begin{document}

\title{Functions on the zeroes of $dx$}

\author{James J. Faran, V}

\address{Department of Mathematics, SUNY at Buffalo, Buffalo, NY, 14260-2900,
USA}

\email{jjfaran@buffalo.edu}

\date{July 23, 2004}

\begin{abstract}
In the model $\mathcal{F}$ of synthetic differential geometry consisting
of sheaves (with respect to open covers) over $\mathbb{F}$, the opposite
category of the category of closed finitely generated $C^{\infty}$-rings,
any morphism from $S$, the zeroes of the ``amazing right adjoint''
of $dx$, to the real line $R$ extends to a morphism from $R$ to
$R$. This shows that the De Rham cohomology of the space $S$ is
the same as the characteristic cohomology of the ideal generated by
$dx$.
\end{abstract}

\keywords{characteristic cohomology, synthetic differential geometry, amazing
right adjoint}

\maketitle

\section{Motivation and Statement of Results}

In two papers (\cite{BG1,BG2}), Bryant and Griffiths developed the
notion of the characteristic cohomology of an exterior differential
system, which has immediate application to integral invariants of
certain parabolic differential equations. One of the open questions
these papers leave is whether there is an analog for characteristic
cohomology to the De Rham theorem for the usual De Rham cohomology.
The question is difficult in part because characteristic cohomology
is defined as an inverse limit. It is not formed directly from looking
at differential forms on spaces. This paper is an attempt to begin
an analysis of this question using synthetic differential geometry.

Certain models of synthetic differential geometry allow differential
forms to be representable, that is, a differential form on a manifold
is a map from that manifold to a certain generalized space. Usually,
the forms under consideration in characteristic cohomology are one-forms
and their derivatives. The ideals are generated (as closed differential
ideals) by one-forms. Representing these one-forms as maps from the
manifold to a generalized space, we can consider the zeroes of these
maps, which together form a generalized subspace of the original manifold.
The long range goal is to show that De Rham's theorem holds on this
generalized space. Here we just provide a simplest case example where
the De Rham cohomology of the generalized subspace is the characteristic
cohomology of the differential ideal.

Now this is a oversimplified picture, since it does not involve the
infinite prolongation necessary to an honest analysis of characteristic
cohomology, but is sufficient in the very elementary case approached
in this paper. Here we consider the ideal of forms on the real line
generated by the standard one form $dx$. The cohomology group of
greatest interest here is that in dimension zero consisting of smooth
functions on the real line whose exterior derivatives are multiples
of $dx$. This clearly is the set of smooth functions on the real
line. Within synthetic differential geometry, we can take the {}``amazing
right adjoint'' of the differential form $dx$, that is, $\check{dx}:R\to\Lambda_{0}$,
and form the subobject $S\subset R$ of the zeroes of $\check{dx}$.
The object of what follows is to prove the following result.

\begin{thm}
In the model $\mathcal{F}$ of synthetic differential geometry consisting
of sheaves (with respect to open covers) over $\mathbb{F}$, the opposite
category of the category of closed finitely generated $C^{\infty}$-rings,
any morphism from $S$, the zeroes of the ``amazing right adjoint''
of $dx$, to the real line $R$ extends to a morphism from $R$ to
$R$.
\end{thm}

\newcommand{\R}{\mathbb{R}}
Since any morphism from $R$ to $R$ corresponds to a smooth map from
$\R$ to $\R$, and since any $R$-valued morphism on $S$ has exterior
derivative zero, this shows that the De Rham cohomology of $S$ is
precisely the characteristic cohomology of the ideal under consideration.
Thus in this simplest case we have found a space that carries the
characteristic cohomology.

The outline of the remainder of this paper is as follows. In section
2, the space $S$ is described as a sheaf. In section 3, we show that
the global section of a morphism $\varphi:S\to R$ gives us a smooth
map $\varphi_{\R}:\R\to\R$ by first defining a derivative morphism
$\varphi':S\to R$ and then showing that the global section of $\varphi'$
is the derivative of the global section of $\varphi$. The global
section $\varphi_{\R}$ gives us our extension of $\varphi$ to $R$,
as is shown in section 4 by showing that any morphism which vanishes
at all the points of $S$ (a point being a map from the terminal space
$*$ to $S$; note that any point of $R$ is a point of $S$) must
be the zero morphism. The arguments in the last two sections are accomplished
by first reducing to questions about certain quotients of power series
rings (this made possible because of the category $\mathcal{F}$ we
are using) and then taking sufficient examples and splitting the singularities
into point-determined spaces. In the last section, the sufficiency
of the set of examples considered is shown using the Briançon-Skoda
theorem.

The basic reference for synthetic differential geometry for this paper
is the book by Moerdijk and Reyes (\cite{MR}). An attempt has been
made to follow their notation as closely as possible. The reader may
need to refer to that text both as background material and for certain
results quoted below.

\thanks{I would like to thank my collegues Bill Lawvere, Steve Schanuel and
Don Schack for their encouragement, support and guidance during this
project, and John D'Angelo and his colleagues at the University of
Illinois for leading me to the Briançon-Skoda theorem.}

\section{The Zeroes of $dx$}

Let $x:R\rightarrow R$ be the identity, $dx:R^{D}\rightarrow R$
its differential, \[
dx\left(\left(d\mapsto x_{0}+dx_{1}\right)\right)=x_{1},\]
 $\check{dx}:R\rightarrow R_{D}$ the right adjoint of $dx$. We let
$S$ be the equalizer \[
S\hookrightarrow R\begin{array}{c}
\check{dx}\\
\rightrightarrows\\
\check{0}\end{array}R_{D}\]
 where $\check{0}:R\rightarrow R_{D}$ is the right adjoint of the
zero differential form $0:R^{D}\rightarrow R$, \[
0\left(\left(d\mapsto x_{0}+dx_{1}\right)\right)=0.\]

What is $S$? As a sheaf, $S$ is a subsheaf of $R$. \begin{eqnarray*}
R\left(\ell A\right) & = & \hom\left(\ell A,\ell C^{\infty}\left(\mathbb{R}\right)\right)\\
 & \simeq & \hom\left(C^{\infty}\left(\mathbb{R}\right),A\right)\\
 & \simeq & A.\end{eqnarray*}
 An element of $R$ is a morphism $\alpha:\ell A\rightarrow\ell C^{\infty}\left(\mathbb{R}\right)$.
This will be an element of $S$ iff it factors through $S$ iff $\check{0}\circ\alpha=\check{dx}\circ\alpha$
iff $0\circ\alpha^{D}=dx\circ\alpha^{D}$ iff the one form $d\alpha:\ell A^{D}\rightarrow R$
is the zero one-form. 

Now $\alpha:\ell A\rightarrow R$ corresponds to an element $a\in A$.
Given $a$, we want to calculate $d\alpha$. 

Now \[
\ell A^{D}=\kappa\left(i\left(\ell A\right)^{D}\right),\]
 where $\kappa$ is the right adjoint of the injection $i:\mathbb{F\hookrightarrow L}$
(\cite{MR}, p. 59). If $A=C^{\infty}\left(\mathbb{R}^{n}\right)/I,$\[
i\left(\ell A\right)^{D}=C^{\infty}\left(\mathbb{R}^{n}\times\mathbb{R}^{n}\right)/\left(I\left(x\right),\sum y_{i}\frac{\partial f}{\partial x_{i}}\mid f\in I\right).\]
 Let \[
I^{D}=\left(I\left(x\right),\sum y_{i}\frac{\partial f}{\partial x_{i}}\mid f\in I\right).\]
 We obtain the reflection $\kappa\left(i\left(\ell A\right)^{D}\right)$
by replacing $I^{D}$ by $\overline{I^{D}}$, the smallest closed
ideal containing $I^{D}$ (\cite{MR}, p. 49). Note that $\overline{I^{D}}$
is not always $I^{D}$. 

Thus, in $\mathbb{F}$, \[
\ell A^{D}=\ell\left(C^{\infty}\left(\mathbb{R}^{n}\times\mathbb{R}^{n}\right)/\overline{I^{D}}\right)\]
 (\emph{cf.} \cite{MR}, p. 71, II.1.16). 

Now let $Y:\mathbb{F}\rightarrow\text{Sets}^{\mathbb{F}_{\text{op}}}$
be the Yoneda embedding, \[
Y\left(\ell A\right)=\hom_{\mathbb{F}}\left(-,\ell A\right).\]

\begin{lem}
$Y\left(\ell A^{D}\right)=Y\left(\ell A\right)^{Y\left(D\right)}$.
\end{lem}
\begin{proof}
From p. 74 of \cite{MR},\begin{eqnarray*}
Y\left(\ell A\right)^{Y\left(D\right)}\left(\ell B\right) & = & \text{Sets}^{\mathbb{F}_{\text{op}}}\left(\ell B\times Y\left(D\right),\ell A\right)\\
 & = & \mathbb{F}\left(\ell B\times D,\ell A\right)\qquad\text{(since }Y\text{ is full)}\\
 & = & \mathbb{F}\left(\ell B,\ell A^{D}\right)\\
 & = & Y\left(\ell A^{D}\right)\left(\ell B\right).\end{eqnarray*}

\end{proof}
So, suppressing the $Y$, $\ell A^{D}=\ell\left(C^{\infty}\left(\mathbb{R}^{n}\times\mathbb{R}^{n}\right)/\overline{I^{D}}\right)$
in $\text{Sets}^{\mathbb{F}_{\text{op}}}$. The resulting presheaf
is a sheaf (see Appendix 1, \S2.7 (c), p. 352 in \cite{MR}), and
so the same result holds in $\mathcal{F}$. 

Let us return to the problem of given $\alpha\in\mathcal{F}\left(\ell A,R\right)$
calculating $\alpha^{D}\in\mathcal{F}\left(\ell A^{D},R^{D}\right)$.
First note $\mathcal{F}\left(\ell A,R\right)=\mathbb{F}\left(\ell A,R\right)$
since the Yoneda embedding is full and faithful, so $\alpha\in\mathbb{F}\left(\ell A,R\right)$. 

If $A=C^{\infty}\left(\mathbb{R}^{n}\right)/I$ ($I$ closed), $\alpha:\ell A\rightarrow R=\ell C^{\infty}\left(\mathbb{R}\right)$
corresponds to $\bar{\alpha}:C^{\infty}\left(\mathbb{R}\right)\rightarrow A$
which corresponds to $a=\bar{\alpha}\left(x\right)\in A$. ($a=\left[f\right]$,
where $f:\mathbb{R}^{n}\rightarrow\mathbb{R}$ is $C^{\infty}$, well-defined
modulo $I$ and $\bar{\alpha}\left(g\right)=\left[f\circ g\right]$.) 

Since $i:\mathbb{F}\rightarrow\mathbb{L}$ is full and faithful, we
can consider $\alpha\in\mathbb{F}\left(\ell A,R\right)$ as $\alpha\in\mathbb{L}\left(i\left(\ell A\right),R\right)$,
calculate $\alpha^{D}\in\mathbb{L}\left(i\left(\ell A\right)^{D},R^{D}\right)$
and then reflect $\kappa\left(\alpha^{D}\right)\in\mathbb{F}\left(\kappa\left(\left(i\left(\ell A\right)\right)^{D}\right),R^{D}\right)=\mathbb{F}\left(\ell A^{D},R^{D}\right).$
Now $i\left(\ell A\right)^{D}$ is isomorphic to $T\left(\ell A\right)=\ell\left(C^{\infty}\left(\mathbb{R}^{n}\times\mathbb{R}^{n}\right)/\left(I\left(x\right),\left\{ \sum y_{i}\frac{\partial f}{\partial x_{i}}\mid f\in I\right\} \right)\right)$,
and the map $\alpha^{D}$ is given as the algebra map $x\mapsto a,y\mapsto\sum y_{i}\frac{\partial a}{\partial x_{i}}$.
The map we want is the composition of this map with the quotient \[
C^{\infty}\left(\mathbb{R}^{n}\times\mathbb{R}^{n}\right)/\left(I\left(x\right),\left\{ \sum y_{i}\frac{\partial f}{\partial x_{i}}\mid f\in I\right\} \right)\rightarrow C^{\infty}\left(\mathbb{R}^{n}\times\mathbb{R}^{n}\right)/\overline{I^{D}}.\]
 Thus $\alpha\in S$ if and only if $\sum y_{i}\frac{\partial f}{\partial x_{i}}$
has Taylor series in the Taylor series of functions in $I^{D}$ everywhere
in $Z\left(I^{D}\right)$. This is equivalent to the following condition.

\begin{condition}
\label{con:Aflat}For all $i=1,\ldots,n$, and at all points of $Z\left(I^{D}\right)$,
$\sum_{i=1}^{n}y^{i}\frac{\partial f}{\partial x_{i}}$ has Taylor
series in the ideal of Taylor series of functions generated by the
Taylor series of functions in $I$ and the Taylor series of the functions
$\sum_{i=1}^{n}y^{i}\frac{\partial g}{\partial x_{i}}$ for $g\in I$.
\end{condition}
\begin{defn}
Let\[
A_{\flat}=\left\{ \left[f\right]\in A:f\textrm{ satisfies Condition \ref{con:Aflat}}\right\} \]
 Note that $A\mapsto A_{\flat}$ is a functor, whose maps are given
by restriction, that is, if $\alpha:A\rightarrow B$, $\alpha_{\flat}:A_{\flat}\rightarrow B_{\flat}$
is simply the restriction of $\alpha$ to $A_{\flat}$. 
\end{defn}
We have shown the following.

\begin{prop}
$S(A)=A_{\flat}.$
\end{prop}
Condition \ref{con:Aflat} as it stands is in need of simplification.
This simplification comes from the simple way in which the variables
$y^{i}$ appear. Suppose $x^{0}\in Z\left(I\right)$. Then $\left(x^{0},0\right)\in Z\left(I^{D}\right)$.
Calculating Taylor series at such points,\[
T_{\left(x^{0},0\right)}\left(\sum_{i=1}^{n}y_{i}\frac{\partial f}{\partial x_{i}}\right)=\sum_{i=1}^{n}y_{i}T_{x^{0}}\left(\frac{\partial f}{\partial x_{i}}\right),\]
and for $g\in I$,\[
T_{\left(x^{0},0\right)}\left(g\right)=T_{x^{0}}\left(g\right),\]
\[
T_{\left(x^{0},0\right)}\left(\sum_{i=1}^{n}y_{i}\frac{\partial g}{\partial x_{i}}\right)=\sum_{i=1}^{n}y_{i}T_{x^{0}}\left(\frac{\partial g}{\partial x_{i}}\right).\]
Thus the condition is that there are power series $b_{j},c_{k}$ and
functions $g_{j},h_{k}\in I$ such that\[
\sum_{i=1}^{n}y_{i}T_{x^{0}}\left(\frac{\partial f}{\partial x_{i}}\right)=\sum_{j}b_{j}T_{x^{0}}\left(g_{j}\right)+\sum_{k}c_{k}\sum_{i=1}^{n}y_{i}T_{x^{0}}\left(\frac{\partial h_{k}}{\partial x_{i}}\right).\]
Equating coefficients of $y_{i}$, we get that\[
T_{x^{0}}\left(\frac{\partial f}{\partial x_{i}}\right)=\sum_{j}b_{ji}T_{x^{0}}\left(g_{j}\right)+\sum_{k}c_{k0}T_{x^{0}}\left(\frac{\partial h_{k}}{\partial x_{i}}\right),\]
where $b_{ji}=b_{ji}\left(x\right),c_{k0}=c_{ko}\left(x\right)$ are
the power series in $x$ which are the coefficients of $y_{i}$ in
$b_{j}$ and $1$ in $c_{k}$ as we expand those series first as series
in $y$. Notice, therefore, that\[
T_{x^{0}}\left(\frac{\partial}{\partial x_{i}}\left(f-\sum_{k}c_{k0}h_{k}\right)\right)=\sum_{j}b_{ji}T_{x^{0}}\left(g_{j}\right)-\sum_{k}\frac{\partial c_{k}}{\partial x_{i}}T_{x^{0}}\left(h_{k}\right).\]
Thus the Condition \ref{con:Aflat} implies the following condition.

\begin{condition}
\label{con:Aflat2}For every point $x^{0}$ of $Z\left(I\right)$,
$a$ has a representative whose partial derivatives have Taylor series
at $x^{0}$ in the ideal of Taylor series of functions in $I$.
\end{condition}
Suppose now that $x^{0}\in Z\left(I\right)$ and $\sum_{i=1}^{n}y_{i}^{0}\frac{\partial g}{\partial x_{i}}\left(x^{0}\right)=0$,
so $\left(x^{0},y^{0}\right)\in Z\left(I^{D}\right)$. Expanding in
a Taylor series about $\left(x^{0},y^{0}\right)$ and using Condition
\ref{con:Aflat2}, it is easy to see that Condition \ref{con:Aflat}
is satisfied. Thus we shall use Condition \ref{con:Aflat2} as our
description of $A_{\flat}$. Note that it is most useful when $A$
is a local $C^{\infty}$-ring (since then there is only one point
$x^{0}\in Z\left(I\right)$), and is more useful as a description
of a special property of elements of $A_{\flat}$ than as a tool for
determining whether a particular element $a\in A$ is an element of
$A_{\flat}$.

\section{Smoothness of the Core}

Let $\phi:S\rightarrow R$ be a morphism. Then since the global sections
of $S$ are precisely the points of $\mathbb{R}$, the global sections
of $\phi$ give a map \[
\phi_{\mathbb{R}}=\Gamma(\phi):\mathbb{R}\rightarrow\mathbb{R},\]
 which will be called the \emph{core} of the morphism $\phi$. The
objective of this section is to show that the core is a smooth map. 

This will be done by defining a morphism $\phi^{\prime}:S\rightarrow R$
and showing that the core of $\phi^{\prime}$ is the derivative of
the core of $\phi$. Since $\phi^{\prime}$ is also a morphism from
$S$ to $R$, we can iterate this argument to show that the core of
$\phi$ is infinitely differentiable. However, before we proceed with
the proof we need some results about $C^{\infty}$-rings.

\subsection{Algebraic Preliminaries}

First, we shall need the left exactness of the coproduct. Alas, this
is not true in generality. It is true, however, if we consider homomorphisms
between Weil algebras.

Let $A=C^{\infty}\left(\mathbb{R}^{n}\right)/I$ be a finitely generated
closed $C^{\infty}$-ring, and let $W=C^{\infty}\left(\mathbb{R}^{m}\right)/J$
be a Weil algebra. We may assume that $J$ is contained in the maximal
ideal of functions vanishing at the origin. (For basic properties
of Weil algebras see \cite{MR}, p. 35\emph{ff}.) Their coproduct
as $C^{\infty}$-rings\[
A\otimes_{\infty}W=C^{\infty}\left(\mathbb{R}^{n}\times\mathbb{R}^{m}\right)/\left(I,J\right).\]
$W$ has a basis $\left\{ w_{1},\ldots,w_{p}\right\} $, and there
are polynomials \[
h_{i}\in\mathbb{R}\left[y_{1},\ldots,y_{m}\right]\subset C^{\infty}\left(\mathbb{R}^{m}\right)\]
 representing the $w_{i}$, $i=1,\ldots,p$.

Let $f\in C^{\infty}\left(\mathbb{R}^{n}\times\mathbb{R}^{m}\right)$.
Expanding $f$ in a Taylor series in the last $m$ variables, and
using the fact that $J$ contains some power of the maximal ideal
of functions vanishing at the origin, we can write\[
f\equiv\sum_{\left|k\right|\le N}f_{jk}\left(x\right)y^{k}\qquad\text{mod }J,\]
where $k$ is a multi-index. We can then write, mod $J$, each $y^{k}$
as a linear combination of the $h_{i}$, so\[
f\equiv\sum_{i=1}^{p}f_{i}\left(x\right)h_{i}\left(y\right)\qquad\text{mod }J.\]

So suppose\[
f=\sum_{i=1}^{p}f_{i}\left(x\right)h_{i}\left(y\right).\]
Then the Taylor series of $f$ at $\left(x_{0},0\right)$,\[
T_{\left(x_{0},0\right)}f\left(x,y\right)=\sum_{i=1}^{p}T_{x_{0}}f_{i}(x)h_{i}\left(y\right).\]
Now we can write any Taylor series uniquely as $\sum_{i=1}^{p}g_{i}\left(x\right)h_{i}\left(y\right)$
modulo $T_{0}J$, just as we did for functions. So $T_{\left(x_{0},0\right)}f\in T_{\left(x_{0},0\right)}\left(I,J\right)$
if and only if $T_{x_{0}}f_{i}\in T_{x_{0}}I$. So, since $I$ is
a closed ideal, if for all $(x_{0},0)\in Z\left(I,J\right)$, $T_{\left(x_{0},0\right)}f\in T_{\left(x_{0},0\right)}\left(I,J\right)$,
then each $f_{i}\in I$, so $f\in\left(I,J\right)$.

Hence $\left(I,J\right)$ is closed and the coproduct of $A$ and
$W$ in the category of finitely generated closed $C^{\infty}$-rings,\[
A\overline{\otimes_{\infty}}W=C^{\infty}\left(\mathbb{R}^{n}\times\mathbb{R}^{m}\right)/\left(I,J\right),\]
which we can identify with\[
\left\{ \sum_{i=1}^{p}f_{i}\left(x\right)h_{i}\left(y\right)\right\} /I.\]

Now suppose $\varphi:W\to W'=C^{\infty}\left(\mathbb{R}^{m'}\right)/J'$
is a homomorphism of Weil algebras. There are bases $\left\{ w_{1},\ldots,w_{p}\right\} $
of $W$ and $\left\{ w_{1}',\ldots,w_{p'}'\right\} $ of $W'$ so
that\[
\varphi\left(w_{i}\right)=\left\{ \begin{array}{cl}
w_{i}' & \text{if }i=1,\ldots,r,\\
0 & \text{if }i=r+1,\ldots,p'.\end{array}\right.\]
Now if $h_{i}'\left(y'\right)$ are polynomials representing $w_{i}'$,
the map\[
\text{id}\overline{\otimes_{\infty}}\varphi:A\overline{\otimes_{\infty}}W\to A\overline{\otimes_{\infty}}W'\]
will map $\left[\sum_{i=1}^{p}f_{i}\left(x\right)h_{i}\left(y\right)\right]$
to $\left[\sum_{i=1}^{r}f_{i}\left(x\right)h_{i}'\left(y'\right)\right]$.
Thus the element $\left[\sum_{i=1}^{p}f_{i}\left(x\right)h_{i}\left(y\right)\right]$
will be in the kernel of $\text{id}\overline{\otimes_{\infty}}\varphi$
if and only if $f_{i}\in I$ for $i=1,\ldots,r$. We have shown the
following.

\begin{prop}
\label{pro:Exactness of Coproduct}Suppose $A$ is a finitely generated
closed $C^{\infty}$-rings, $W$ and $W'$ are Weil algebras, and
$\varphi:W\rightarrow W'$ is a $C^{\infty}$-ring homomorphism. Then\[
\text{id}\overline{\otimes_{\infty}}\varphi:A\overline{\otimes_{\infty}}W\rightarrow A\overline{\otimes_{\infty}}W'\]
 has kernel the ideal generated by $\ker\varphi\overline{\otimes_{\infty}}1$.
\end{prop}
The next result we need is an injectivity result, this time not involving
the coproduct, but rather to a certain pushout. As a preliminary,
we show the following.

\begin{lem}
Let $A$ and $B$ be fintely generated closed $C^{\infty}$-rings.
Let $a\in A$, $a\ne0$ and $b\in B$, $b\ne0$. Then $a\overline{\otimes_{\infty}}b\ne0$
in $A\overline{\otimes_{\infty}}B$.
\end{lem}
\begin{proof}
Suppose $A=C^{\infty}\left(\mathbb{R}^{n}\right)/I$ and $B=C^{\infty}\left(\mathbb{R}^{m}\right)/J$,
so $A\overline{\otimes_{\infty}}B=C^{\infty}\left(\mathbb{R}^{n}\times\mathbb{R}^{m}\right)/\overline{\left(I,J\right)}$.
It suffices to show that there is some $\left(p,q\right)\in Z\left(\overline{\left(I,J\right)}\right)$,
the set of common zeros of the ideal $\overline{\left(I,J\right)}$,
for which the Taylor series of a representative of $a\overline{\otimes_{\infty}}b$
at $\left(p,q\right)$ is not in the ideal of Taylor series consisting
of Taylor series of elements of $\overline{\left(I,J\right)}$ at
$\left(p,q\right)$. However, since $a\ne0$, there is some $p\in\mathbb{R}^{n}$
at which the Taylor series of a representative of $a$ does not lie
in $T_{p}I$, the ideal consisting of Taylor series of elements of
$I$ at $p$, and since $b\ne0$, there is some $q\in\mathbb{R}^{m}$
at which the Taylor series of a representative of $b$ does not lie
in $T_{q}J$.

Thus we may assume that\[
A=\mathbb{R}\left[\left[x_{1},\ldots,x_{n}\right]\right]/I,\]
 \[
B=\mathbb{R}\left[\left[y_{1},\ldots,y_{m}\right]\right]/J,\]
 and hence\[
A\overline{\otimes_{\infty}}B=\mathbb{R}\left[\left[x_{1},\ldots,x_{n},y_{1},\ldots,y_{m}\right]\right]/\left(I,J\right).\]
By the Krull intersection theorem (see, for example, \cite{L}, p.
430, Corollary 5.7), there is some $N_{1}\in\mathbb{N}$ so that $a\ne0$
in $A/\mathfrak{m}_{A}^{N_{1}}$, where $\mathfrak{m}_{A}$ is the
maximal ideal of $A$. Similarly, there is some $N_{2}\in\mathbb{N}$
so that $b\ne0$ in $B/\mathfrak{m}_{B}^{N_{2}}$. Therefore it suffices
to show that $a\overline{\otimes_{\infty}}b\ne0$ in $A/\mathfrak{m}_{A}^{N_{1}}\overline{\otimes_{\infty}}B/\mathfrak{m}_{B}^{N_{2}}$.
But $A/\mathfrak{m}_{A}^{N_{1}}$ and $B/\mathfrak{m}_{B}^{N_{2}}$
are both Weil algebras, for which $\overline{\otimes_{\infty}}$ is
simply the standard tensor product, and it is clear that $a\ne0$
and $b\ne0$ implies $a\otimes b\ne0$.
\end{proof}
\begin{prop}
\label{pro:pushout injects}Let $A=C^{\infty}\left(\mathbb{R}^{n}\right)/I$
and $B=C^{\infty}\left(\mathbb{R}^{m}\right)/J$ be two finitely generated
closed $C^{\infty}$-rings, and for some $n\in\mathbb{N}$, let $a_{0}\in A$,
$b_{0}\in B$ be two elements with $a_{0}^{n}=0$, $a_{0}^{n-1}\ne0$
, $b_{0}^{n}=0$, $b_{0}^{n-1}\ne0$. Corresponding to these two elements
are maps $C^{\infty}\left(\mathbb{R}\right)\rightarrow A$ and $C^{\infty}\left(\mathbb{R}\right)\rightarrow B$
sending the standard generator of $C^{\infty}\left(\mathbb{R}\right)$
to $a_{0}$ and $b_{0}$. The pushout of these two maps we write as
$A\overline{\otimes_{C^{\infty}\left(\mathbb{R}\right)}}B$ and is
seen to be $A\overline{\otimes_{\infty}}B/\left(a_{0}-b_{0}\right)$.
(Here we have identified $a_{0}\in A$ and its image in the coproduct
$A\overline{\otimes_{\infty}}B$ under the canonical map $A\rightarrow A\overline{\otimes_{\infty}}B$,
and similarly for $b_{0}$.) 

The canonical map $A\rightarrow A\overline{\otimes_{C^{\infty}\left(\mathbb{R}\right)}}B$
is then an injection.
\end{prop}
\begin{proof}
Let $a\in A$ have image zero in $A\overline{\otimes_{C^{\infty}\left(\mathbb{R}\right)}}B$.
We start by noting that since $B$ is a closed $C^{\infty}$-ring,
it has a {}``point'', that is, there is a homomorphism $B\to\mathbb{R}$.
Since $b_{0}^{n}=0$, any such point is also a point of $B/\left(b_{0}\right)$,
so there is a homomorphism $p:B/\left(b_{0}\right)\to\mathbb{R}$.
Thus there is a map $1\overline{\otimes_{\infty}}p:A\overline{\otimes_{\infty}}B/\left(b_{0}\right)\to A\overline{\otimes_{\infty}}\mathbb{R}=A$.

Since $a$ has image zero in $A\overline{\otimes_{C^{\infty}\left(\mathbb{R}\right)}}B$,\[
a\overline{\otimes_{\infty}}1=m\left(a_{0}\overline{\otimes_{\infty}}1-1\overline{\otimes_{\infty}}b_{0}\right).\]
Then in $A\overline{\otimes_{\infty}}B/\left(b_{0}\right)$,\[
a\overline{\otimes_{\infty}}1=m\left(a_{0}\overline{\otimes_{\infty}}1\right).\]
Applying the map $1\overline{\otimes_{\infty}}p$,\[
a=\left(1\overline{\otimes_{\infty}}p\right)\left(q\left(m\right)\right)\cdot a_{0},\]
where $q:A\overline{\otimes_{\infty}}B\to A\overline{\otimes_{\infty}}B/\left(b_{0}\right)$
is the quotient map. Therefore\[
a\overline{\otimes_{\infty}}1=\left(\left(1\overline{\otimes_{\infty}}p\right)\left(q\left(m\right)\right)\overline{\otimes_{\infty}}1\right)\cdot\left(a_{0}\overline{\otimes_{\infty}}1\right)\]
which we will write as\[
a\overline{\otimes_{\infty}}1=\left(a'\overline{\otimes_{\infty}}1\right)\cdot\left(a_{0}\overline{\otimes_{\infty}}1\right).\]
Hence \[
0=\left(m-a'\overline{\otimes_{\infty}}1\right)\cdot\left(a_{0}\overline{\otimes_{\infty}}1\right)-m\left(1\overline{\otimes_{\infty}}b_{0}\right).\]
Writing $m'=m-a'\overline{\otimes_{\infty}}1$,\begin{equation}
0=m'\left(a_{0}\overline{\otimes_{\infty}}1\right)-m'\left(1\overline{\otimes_{\infty}}b_{0}\right)+a'\overline{\otimes_{\infty}}b_{0}.\label{eq: vanishing 1}\end{equation}
Multiplying this by $\left(a_{0}\overline{\otimes_{\infty}}1\right)^{n-1}$,\[
0=m'\left(a_{0}^{n}\overline{\otimes_{\infty}}1\right)-m'\left(a_{0}^{n-1}\overline{\otimes_{\infty}}b_{0}\right)+a'a_{0}^{n-1}\overline{\otimes_{\infty}}b_{0}\]
\begin{eqnarray*}
0 & = & m'\left(a_{0}^{n}\overline{\otimes_{\infty}}1\right)-m'\left(a_{0}^{n-1}\overline{\otimes_{\infty}}b_{0}\right)+a'a_{0}^{n-1}\overline{\otimes_{\infty}}b_{0}\\
 & = & -m'\left(a_{0}^{n-1}\overline{\otimes_{\infty}}b_{0}\right)+a'a_{0}^{n-1}\overline{\otimes_{\infty}}b_{0}\end{eqnarray*}
since $a_{0}^{n}=0$. We can now use (\ref{eq: vanishing 1}) to replace
$m'\left(a_{0}\overline{\otimes_{\infty}}1\right)$ by the expression
$m'\left(1\overline{\otimes_{\infty}}b_{0}\right)-a'\overline{\otimes_{\infty}}b_{0}$,
obtaining\begin{eqnarray*}
0 & = & -\left(m'\left(1\overline{\otimes_{\infty}}b_{0}\right)-a'\overline{\otimes_{\infty}}b_{0}\right)\left(a_{0}^{n-2}\overline{\otimes_{\infty}}b_{0}\right)+a'a_{0}^{n-1}\overline{\otimes_{\infty}}b_{0}\\
 & = & -m'a_{0}^{n-2}\overline{\otimes_{\infty}}b_{0}^{2}+a'a_{0}^{n-2}\overline{\otimes_{\infty}}b_{0}^{2}+a'a_{0}^{n-1}\overline{\otimes_{\infty}}b_{0}.\end{eqnarray*}
Repeating this substitution the appropriate number of times,\[
0=-m'\left(1\overline{\otimes_{\infty}}b_{0}^{n}\right)+a'\overline{\otimes_{\infty}}b_{0}^{n}+\sum_{k=1}^{n-1}a'a_{0}^{n-k}b_{0}^{k},\]
hence, since $b_{0}^{n}=0$,\begin{equation}
0=\sum_{k=1}^{n-1}a'a_{0}^{n-k}\overline{\otimes_{\infty}}b_{0}^{k}.\label{eq: vanishing2}\end{equation}
Multiplying this by $a_{0}^{n-2}\overline{\otimes_{\infty}}1$,\[
0=\sum_{k=1}^{n-1}a'a_{0}^{2n-k-2}\overline{\otimes_{\infty}}b_{0}^{k}=a'a_{0}^{n-1}\overline{\otimes_{\infty}}b_{0}^{n-1},\]
since $a_{0}^{n}=0$. Since $b_{0}^{n-1}\ne0$, by the previous lemma
$a'a_{0}^{n-1}=0.$ Using this, we multiply (\ref{eq: vanishing2})
by $a_{0}^{n-2}\overline{\otimes_{\infty}}1$ and similarly obtain
$a'a_{0}^{n-2}=0.$ Repeating this process an appropriate number of
times (multiplying (\ref{eq: vanishing2}) by smaller and smaller
powers of $a_{0}\overline{\otimes_{\infty}1}$), we eventually obtain
$a'a_{0}=0$. However, $a'a_{0}=a$, so $a=0$ as desired.
\end{proof}

\subsection{Existence of the Derivative Morphism}

Recall that the space of quadratic nilpotents $D=\left\{ d\in R:d^{2}=0\right\} $.

\begin{prop}
Given a morphism $\phi:S\to R$, for all $s\in S$ and all $d\in S\cap D$,
\[
\phi(s+d)=\phi(s)+\phi^{\prime}(s)d\]
 for some $\phi^{\prime}(s)\in R$.
\end{prop}
\begin{proof}
To prove this, we start by seeing what it means. $s\in S$ corresponds
to a morphism $\ell A\to S$ which corresponds to $a\in A_{\flat}$.
$d\in S\cap D$ corresponds to a morphism $\ell B\to S\cap D$ which
corresponds to $b\in B_{\flat}$ such that $b^{2}=0$. Then $s+d$
corresponds to $a+b=a\overline{\otimes_{\infty}}1+1\overline{\otimes_{\infty}}b\in\left(A\overline{\otimes_{\infty}}B\right)_{\flat}$.
We need to show that \[
\phi_{A\overline{\otimes_{\infty}}B}(a+b)=\phi_{A}(a)\overline{\otimes_{\infty}}1+\phi^{\prime}(a)\overline{\otimes_{\infty}}b\]
 for some $\phi^{\prime}(a)\in A$. 

We start by considering the quotient \[
A\overline{\otimes_{\infty}}B\to A\overline{\otimes_{\infty}}B/(1\overline{\otimes_{\infty}}b).\]
 From the functoriality of $\phi$, chasing $a+b$ around the commutative
square, we have \[
\left[\phi_{A\overline{\otimes_{\infty}}B}(a+b)\right]=\phi_{A\overline{\otimes_{\infty}}B/(1\overline{\otimes_{\infty}}b)}\left(\left[a\overline{\otimes_{\infty}}1\right]\right).\]
 Now we also have the map $A\to A\overline{\otimes_{\infty}}B\to A\overline{\otimes_{\infty}}B/(1\overline{\otimes_{\infty}}b)$ 
 that sends $a\mapsto a\overline{\otimes_{\infty}}1\mapsto\left[a\overline{\otimes_{\infty}}1\right]$,
which yields \[
\phi_{A\overline{\otimes_{\infty}}B/(1\overline{\otimes_{\infty}}b)}\left(\left[a\overline{\otimes_{\infty}}1\right]\right)=\left[\phi_{A}(a)\overline{\otimes_{\infty}}1\right].\]
 Therefore \[
\left[\phi_{A\overline{\otimes_{\infty}}B}(a+b)\right]=\left[\phi_{A}(a)\overline{\otimes_{\infty}}1\right],\]
 so \[
\phi_{A\overline{\otimes_{\infty}}B}(a+b)=\phi_{A}(a)\overline{\otimes_{\infty}}1+\psi_{A,B}(a,b)(1\overline{\otimes_{\infty}}b)\]
 for some $\psi_{A,B}(a,b)$. We need to show that not only can $\psi_{A,B}(a,b)$
be chosen independent of $b$ and $B$, but is of the form $\phi^{\prime}(a)\overline{\otimes_{\infty}}1$.

To this end, suppose we have some element $b'$ of a closed $C^{\infty}$-ring $B'$
 with $b'^{2}=0$. We will take the following as a fundamental example.

Let \[
B_{0}=C^{\infty}\left(\mathbb{R}^{2}\right)/(5x^{4}+3x^{2}y^{3},3x^{3}y^{2}+5y^{5},y^{6}).\]
 Then $b_{0}=x^{5}+x^{3}y^{3}+y^{5}$ is non-zero in $B_{0}$, but
$b_{0}^{2}=0$ in $B_{0}$. Moreover, $b_{0}\in\left(B_{0}\right)_{\flat}$.

Now consider the maps from $A\otimes_{\infty}B$ and $A\otimes_{\infty}B_{0}$
to $C=A\otimes_{\infty}B\otimes_{\infty}B_{0}/(b-b_{0})$. As seen
above, these maps are injections. Both $a+b$ and $a+b_{0}$ get sent
to the same element of $C$. By the functoriality of $\phi$ we have
\[
\left[\phi(a)+\psi_{A,B}(a,b)b\right]=\left[\phi(a)+\psi_{A,B_{0}}(a,b_{0})b\right]\]
 in $C$. However, the kernel of the map $B_{0}\to B_{0}/(b_{0})$,
the ideal $\left(b_{0}\right)$, consists entirely of scalar multiples
of $b_{0}$: $(b_{0})=\mathbb{R}b_{0}$. By the exactness of the coproduct
(Proposition \ref{pro:Exactness of Coproduct}), the kernel of the
map $A\overline{\otimes_{\infty}}B_{0}$ to $A\overline{\otimes_{\infty}}B_{0}/\left(b_{0}\right)$
is $A\overline{\otimes_{\infty}}b_{0}$. Thus $\psi_{A,B_{0}}(a,b_{0})b=\phi_{A}^{\prime}(a)b$,
where $\phi_{A}^{\prime}(a)\in A$.

We then have \[
\left[\psi_{A,B}(a,b)b\right]=\left[\phi_{A}^{\prime}(a)b\right]\]
 in $C$, and so, by the injectivity of the pushout (Proposition \ref{pro:pushout injects}),
\[
\psi_{A,B}(a,b)b=\phi_{A}^{\prime}(a)b\]
 in $A$. Since $\phi_{A}^{\prime}(a)\in A$ is unique, this defines
a collection of maps. Moreover, since the construction is natural,
it defines a morphism $\phi^{\prime}:S\to R$.
\end{proof}

\subsection{The Derivative of the Core}

Our goal now is to show that the derivative of the core (as a function
from $\mathbb{R}$ to $\mathbb{R}$) is the core of the derivative
(the morphism just defined). This is done by taking a simple locus,
where we can see clearly how the derivative must operate, and {}``splitting''
the singularity, obtaining loci consisting of points. Since these
{}``split'' loci consist only of points, the behavior of the original
morphism $\phi$ is governed completely by the behavior of its core.
Because the splitting is done smoothly, we can relate difference quotients
to the action of the derivative.

\subsubsection{The Simple Locus}

We consider a quotient of the ring of real power series in two variables.
Let \[
A=\mathbb{R}[[x,y]]/\left(5x^{4}+2xy^{2},2x^{2}y+5y^{4}\right).\]
 Note that the two functions here are the two partial derivatives
of $f(x,y)=x^{5}+x^{2}y^{2}+y^{5}$, so $[f]\in A_{\flat}$. Moreover,
$[f]\ne0$ in $A$. $A$ has dimension $11$ as a real vector space,
and the ideal of $A$ generated by $[f]$ consists entirely of real
multiples of $[f]$. 

By considering the quotient map $A\rightarrow A/([f])$, we see that
there is a function $\lambda:\mathbb{R}\rightarrow\mathbb{R}$ such
that for all $c\in\mathbb{R}\subset A$, \[
\phi(c+[f])=\phi(c)+\lambda(c)[f].\]
 Indeed, $\lambda(c)=\phi^{\prime}(c)$, the value of the derivative
morphism at $c\in A$. Thus, we wish to show that we obtain the same
$\lambda(c)$ by taking the usual difference quotient.

\subsubsection{The Splitting}

To split the singularity, we need to move away from real power series.
Note the the common real zeros of $5x^{4}+2xy^{2}$ and $2x^{2}y+5y^{4}$
are $(0,0)$ and $(-2/5,-2/5)$. Thus $C^{\infty}\left(\mathbb{R}^{2}\right)/I$,
where $I$ is the ideal generated by $5x^{4}+2xy^{2}$, $2x^{2}y+5y^{4}$
and all functions vanishing on the disk of radius $1/2$ around the
origin, is precisely our ring $A$.

We perturb the defining equations in the following way. Let\[
g_{1}=x\left[5\left(x^{3}+ax^{2}-a^{4}x-a^{5}-\frac{2}{5}a^{4}\right)+2y^{2}\right],\]
\[
g_{2}=y\left[5\left(y^{3}+ay^{2}-a^{4}y-a^{5}-\frac{2}{5}a^{4}\right)+2x^{2}\right].\]
When the parameter $a=0$, $g_{1}$ and $g_{2}$ are the two important
defining functions of $A$. We want to show that in a neighborhood
of the origin for $a\ne0$ small, the space defined by the vanishing
of $g_{1}$ and $g_{2}$ is the union of $11$ points.

\begin{lem}
In a sufficiently small neighborhood of the origin, for all $a>0$
sufficiently close to $0$, $g_{1}$ and $g_{2}$ have precisely $11$
common zeroes, and $dg_{1}\wedge dg_{2}\ne0$ at those points.
\end{lem}
\begin{proof}
First, let us look at the zeroes along the $y$-axis, $\left\{ x=0\right\} $.
Along this axis $g_{1}$ vanishes, so we only have to consider the
vanishing of $g_{2}$. However, along the $y$-axis, $g_{2}=5yh\left(y\right)$,
where\[
h\left(y\right)=y^{3}+ay^{2}-a^{4}y-a^{5}-\frac{2}{5}a^{4}.\]
It is easy to see, though, that for small $a$, \begin{eqnarray*}
h(-2a) & = & -4a^{3}-\frac{2}{5}a^{4}+a^{5}<0,\\
h\left(-\frac{1}{2}a\right) & = & \frac{1}{8}a^{3}-\frac{2}{5}a^{4}-\frac{1}{2}a^{5}>0,\\
h\left(0\right) & = & -\frac{2}{5}a^{4}-a^{5}<0,\\
h\left(a\right) & = & 2a^{3}-\frac{2}{5}a^{4}-2a^{5}>0.\end{eqnarray*}
Thus $h$ has three zeroes near the origin, which gives us three zeroes
of the form $\left(0,y\right)$ with $y\ne0$. By symmetry, there
are also three zeroes of the form $\left(x,0\right)$ with $x\ne0$,
as well as the zero at the origin. Noting that along the $y$-axis,
$dg_{1}$ is a non-zero multiple of $dx$, and since the zeroes of
$yh\left(y\right)$ are discrete, and applying symmetry, at each of
these seven common zeroes, $dg_{1}\wedge dg_{2}\ne0$.

Next, note there are four zeroes at the points $\left(\pm a^{2},\pm a^{2}\right)$,
giving eleven zeroes near the origin.

Finally, looking at all the zeroes over the complex domain, when $a=0$
there are five zeroes at $\left(-\frac{2}{5}\varepsilon^{2},-\frac{2}{5}\varepsilon^{3}\right)$,
where $\varepsilon$ is a fifth root of unity. These are discrete
zeroes, and so for small $a$ there still are five discrete (complex)
zeroes near these points. This makes sixteen zeroes in all, precisely
the number possible by Bezout's theorem. Since these last five are
away from the origin, we are left with the first eleven near the origin.
\end{proof}
In the sequel, we will need more detailed information about the zeroes
of $h(y)=y^{3}+ay^{2}-a^{4}y-a^{5}-\frac{2}{5}a^{4}.$ Essentially,
the zeroes are at $-a$ and $\pm\sqrt{2/5}a^{3/2}$ asympotically
as $a\to0$. In particular, we shall need the following.

\begin{lem}
\label{lem:distance between zeroes of h}The zeroes of $yh(y)$ differ
by more than a fixed constant times $a^{3/2}$ as $a$ approaches
$0$.
\end{lem}
\begin{proof}
One way to see this is to evaluate $h$ at $-a\textrm{ and }-a+\frac{3}{5}a^{2}$,
then at $-\frac{5}{8}a^{3/2}\textrm{ and }-\frac{6}{8}a^{3/2}$, and
finally at $\frac{5}{8}a^{3/2}\textrm{ and }\frac{6}{8}a^{3/2}$,
each time expanding the result in powers of $a$. The intermediate
value theorem places a zero in each of the three intervals, and the
lemma follows.
\end{proof}

\subsubsection{The Reasoning}

We wish to show\[
\lambda\left(c\right)=\lim_{\Delta c\to0}\frac{\varphi\left(c+\Delta c\right)-\varphi\left(c\right)}{\Delta c}\]
where $c,\Delta c\in\mathbb{R}$. First, note that it suffices to
consider sequential limits, that is, to show\[
\forall\Delta c_{n}\to0,\lambda\left(c\right)=\lim_{n\to\infty}\frac{\varphi\left(c+\Delta c_{n}\right)-\varphi\left(c\right)}{\Delta c_{n}}.\]
Second, note that it suffices to consider monotone sequences of constant
sign, that is, to show\[
\forall\Delta c_{n}\searrow0,\lambda\left(c\right)=\lim_{n\to\infty}\frac{\varphi\left(c+\Delta c_{n}\right)-\varphi\left(c\right)}{\Delta c_{n}}\]
and\[
\forall\Delta c_{n}\nearrow0,\lambda\left(c\right)=\lim_{n\to\infty}\frac{\varphi\left(c+\Delta c_{n}\right)-\varphi\left(c\right)}{\Delta c_{n}}.\]
Thirdly, it suffices to consider sequences that converge to $0$ infinitely
fast,\[
n^{p}\Delta c_{n}\to0\qquad\forall p\in\mathbb{N},\]
for example by making sure $\left|\Delta c_{n}\right|<\frac{1}{n^{n}}$.

It follows that if we have such a sequence $\Delta c_{n}$ we can
find a $C^{\infty}$ monotone function $a:\mathbb{R}\to\mathbb{R}$
such that $a\left(\frac{1}{n}\right)=\left|\Delta c_{n}\right|^{1/8}$.
$a$ will be infinitely flat at $0\in\mathbb{R}$.

Supposing we have such a sequence $\Delta c_{n}$ and function $a$,
we take $g_{1}$ and $g_{2}$ as before, now thought of as functions
on $\mathbb{R}^{3}=\left\{ \left(x,y,t\right)\right\} $ taking $a=a\left(t\right)$,
or rather on $U\times\mathbb{R}$, where $U$ is a small neighborhood
of the origin in $\mathbb{R}^{2}$, and consider the locus which is
the intersection of $\left\{ g_{1}=g_{2}=0\right\} $ with $\left\{ \left(x,y,t\right):t=0\text{ or }t=\frac{1}{n},n\in\mathbb{N}\right\} $.
The next proposition gives a way of representing elements of the corresponding
$C^{\infty}$-ring, which we shall call $B$.

\begin{prop}
For all $N\in\mathbb{N}$ and for all $C^{\infty}$ functions $f\left(x,y,t\right)$
defined in a sufficiently small neighborhood of $\left(0,0,0\right)$,
there are smooth functions $\lambda_{nm}\left(t\right)$ such that\[
f\left(x,y,t\right)\equiv\sum_{n,m=0}^{2}\lambda_{nm}x^{n}y^{m}+\lambda_{30}x^{3}+\lambda_{03}y^{3}\quad\mod\left(g_{1},g_{2},a^{N}\right).\]
Moreover, if $f\left(x,y,t\right)=f\left(y,x,t\right)$, we can find
such $\lambda_{nm}$ such that $\lambda_{nm}=\lambda_{mn}$.
\end{prop}
\begin{proof}
First note that given any function $f\left(x,y,t\right)$ we can find
$\lambda_{nm}\left(t\right)$ such that\[
f\left(x,y,t\right)\equiv\sum_{n,m=0}^{2}\lambda_{nm}x^{n}y^{m}+\lambda_{30}x^{3}+\lambda_{03}y^{3}\quad\mod\left(\left.g_{1}\right|_{t=0},\left.g_{2}\right|_{t=0}\right).\]
To see this, first note that\begin{eqnarray*}
\left.g_{1}\right|_{t=0} & = & 5x^{4}+2xy^{2},\\
\left.g_{2}\right|_{t=0} & = & 5y^{4}+2x^{2}y.\end{eqnarray*}
In a sufficiently small neighborhood of the origin $\left(0,0\right)$,
both $x^{6}$ and $y^{6}$ are in the ideal generated by these two
functions. If we expand $f$ as a truncated power series (with remainder)
in $x$ and $y$ with coefficients functions of $t$, modulo these
two functions we can eliminate all sufficiently high order terms.
It is then easy to see that the remaining terms can be replaced with
lower order terms, that the new coefficients still depend smoothly
on $t$, and that if $f$ is symmetric, the new representation of
$f$ is also symmetric.

If we do the same replacements, but using $g_{1}$ and $g_{2}$ instead
of $\left.g_{1}\right|_{t=0}$ and $\left.g_{2}\right|_{t=0}$, we
will get\[
f\left(x,y,t\right)\equiv\sum_{n,m=0}^{2}\lambda_{nm}x^{n}y^{m}+\lambda_{30}x^{3}+\lambda_{03}y^{3}+aF_{1}\left(x,y,t\right).\]
$\mod\left(g_{1},g_{2}\right)$. Applying the same procedure now to
$F_{1}$, we get\[
f\left(x,y,t\right)\equiv\sum_{n,m=0}^{2}\lambda_{nm}x^{n}y^{m}+\lambda_{30}x^{3}+\lambda_{03}y^{3}+a^{2}F_{2}\left(x,y,t\right)\]
$\mod\left(g_{1},g_{2}\right)$. The result follows in this fashion
by induction on $N$.
\end{proof}
Now if $\Delta c_{n}>0$, let $f_{0}\left(x,y,t\right)=x^{2}y^{2}$.
If $\Delta c_{n}<0$, let $f_{0}\left(x,y,t\right)=-x^{2}y^{2}$.
In both cases, $f_{0}$ is equivalent modulo $\left.g_{1}\right|_{t=0}$
and $\left.g_{2}\right|_{t=0}$ to a real multiple of $x^{5}+x^{2}y^{2}+y^{5}$,
and so is in $A_{\flat}$. Since it is independent of $t$, the the
remainder of the locus consists of points approaching the origin infinitely
fast, it is in $B_{\flat}$. We therefore can apply our morphism.

We shall describe the case $\Delta c_{n}>0$. The other case is similar.

Let\[
\varphi\left(c+f_{0}\right)=\sum_{n,m=0}^{2}\lambda_{nm}x^{n}y^{m}+\lambda_{30}x^{3}+\lambda_{03}y^{3}+O\left(a^{N}\right).\]
Since the function $c+f_{0}$ is symmetric, we may assume the $\lambda_{nm}$
are symmetric. If we consider the value of $\varphi\left(c+f_{0}\right)$
at the various points of our locus,\begin{eqnarray*}
\varphi\left(c+f_{0}\right)\left(0,0,\frac{1}{n}\right) & = & \varphi\left(c\right)=\lambda_{00},\\
\varphi\left(c+f_{0}\right)\left(x,0,\frac{1}{n}\right) & = & \varphi\left(c\right)=\lambda_{00},\\
\varphi\left(c+f_{0}\right)\left(0,y,\frac{1}{n}\right) & = & \varphi\left(c\right)=\lambda_{00}.\end{eqnarray*}
Thus, for $y$ any of the zeroes of $yh(y)=y\left(y^{3}+ay^{2}-a^{4}y-a^{5}-\frac{2}{5}a^{4}\right)$,\[
\varphi\left(c\right)=\varphi\left(c\right)+\lambda_{01}y+\lambda_{02}y^{2}+\lambda_{03}y^{3}+O\left(a^{N}\right).\]
Since the zeroes of $yh(y)$ all differ by a fixed constant times
$a^{3/2}$, it follows that we can solve for $\lambda_{01},\lambda_{02},\lambda_{03}$
and obtain that they vanish to some high power of $a$. By symmetry,
so do $\lambda_{10},\lambda_{20},\lambda_{30}$.

By symmetry, the other four points give three equations, corresponding
to $\left(a^{2},a^{2}\right),\left(-a^{2},-a^{2}\right)\textrm{ and }\left(a^{2},-a^{2}\right)$:\begin{eqnarray*}
\varphi\left(a^{8}\right) & = & \varphi\left(c\right)+\lambda_{11}a^{4}+\lambda_{12}a^{6}+\lambda_{22}a^{8}+O\left(a^{N'}\right),\\
\varphi\left(a^{8}\right) & = & \varphi\left(c\right)+\lambda_{11}a^{4}-\lambda_{12}a^{6}+\lambda_{22}a^{8}+O\left(a^{N'}\right),\\
\varphi\left(a^{8}\right) & = & \varphi\left(c\right)-\lambda_{11}a^{4}+\lambda_{12}a^{6}+\lambda_{22}a^{8}+O\left(a^{N'}\right).\end{eqnarray*}
Subtracting the first two tells us that $\lambda_{12}$ vanishes to
high order. Then adding the first and third gives us\[
\frac{\varphi\left(a^{8}\right)-\varphi\left(c\right)}{a^{8}}=\lambda_{22}+O\left(a^{N''}\right).\]
So, in particular, since $a\left(\frac{1}{n}\right)^{8}=\Delta c_{n},$\[
\left.\lambda_{22}\right|_{t=0}=\lim_{n\to\infty}\frac{\varphi\left(c+\Delta c_{n}\right)-\varphi\left(c\right)}{\Delta c_{n}}.\]
However, $\left.\lambda_{22}\right|_{t=0}=\varphi'\left(c\right).$
We have shown the following.

\begin{thm}
The derivative of the core is the core of the derivative:\[
\left(\varphi'\right)_{\mathbb{R}}=\left(\varphi_{\mathbb{R}}\right)'.\]

\end{thm}
\begin{cor}
The core $\varphi_{\mathbb{R}}$ is infinitely differentiable.
\end{cor}

\section{The Core Defines the Morphism}

The core of our morphism, $\varphi_{\mathbb{R}}:\mathbb{R}\to\mathbb{R}$,
being smooth, gives rise to a morphism $\varphi_{\mathbb{R}}:R\to R$,
and by composition, a morphism $\varphi_{\mathbb{R}}:S\to R$. To
show that $\varphi$ extends to $R$ it suffices to show that $\varphi=\varphi_{\mathbb{R}}$,
or, equivalently, that $\psi=\varphi-\varphi_{\mathbb{R}}=0$. What
we know about $\psi$ is that its core is the $0$ map, that is, $\psi$
vanishes at all points of $S$. 

Thus it suffices to show that any morphism $\psi:S\to R$ which vanishes
at all points is the zero morphism. In particular, for any locus $\ell A$,
we want to show that $\psi_{\ell A}:S_{\ell A}=A_{\flat}\to R_{\ell A}=A$
is the zero map. Since we are working only with closed rings, to show
that $\psi_{\ell A}\left(a\right)=0$, we just need to show that the
Taylor series at every point of $\ell A$ of a representative of $\psi_{\ell A}\left(a\right)$
vanishes. Therefore, it suffices to consider $C^{\infty}$-rings of
the form $A=\mathbb{R}\left[\left[x_{1},\ldots,x_{n}\right]\right]/I$.

We shall now proceed as follows. First, we apply a version of the
Briançon-Skoda theorem to show that any $a\in A_{\flat}$ which vanishes
at the point of $A$ is nilpotent. Second, we produce a family of
examples of nilpotent flat elements with arbitrary degree of nilpotency.
Third, we provide a splitting of these examples (as we did in the
previous section) into $C^{\infty}$-rings whose elements are determined
by their values at points. Finally, these splittings are used to show
that the image of an arbitrary $a\in A_{\flat}$ under $\psi$ must
be zero.

\subsection{The Briançon-Skoda Theorem}

In \cite{BS}, Briançon and Skoda showed that any germ of a holomorphic
function of several variables is integral over the ideal generated
by its partial derivatives. Their result has been generalized to other
rings. (See, for example, \cite{LS} and \cite{LT}. Many other references
are available in \cite{H}.) We shall need this result for the ring
of real formal power series in $n$ variables. Rather than force the
reader to the great generality of the cited papers, a proof in this
simple case is provided here.

Let $f\in\mathbb{R}\left[\left[x_{1},x_{2},\ldots,x_{n}\right]\right]$
be a real power series with power series partial derivatives $\frac{\partial f}{\partial x_{i}}$.
Let $\Delta(f)=\left(\frac{\partial f}{\partial x_{1}},\ldots,\frac{\partial f}{\partial x_{n}}\right)$
be the ideal generated by the partial derivatives of $f$, \emph{i.e.},
the \textbf{Jacobian ideal} of $f$. The purpose of this note is to
provide a proof of the following result.

\begin{thm}
\label{thm:nilpotency}If $f\in\mathbb{R}\left[\left[x_{1},x_{2},\ldots,x_{n}\right]\right]$
is a real power series with Jacobian ideal $\Delta(f)=\left(\frac{\partial f}{\partial x_{1}},\ldots,\frac{\partial f}{\partial x_{n}}\right)$
and if the constant term of $f$ vanishes (that is, $f(0)=0$), then
$f$ is nilpotent in $\mathbb{R}\left[\left[x_{1},x_{2},\ldots,x_{n}\right]\right]/\Delta(f)$,
in other words $f$ is in the radical of $\Delta(f)$.
\end{thm}
To prove this, first recall the following theorem (Corollary 2.3 of
Chapter X, \S 2 of \cite{L}).

\begin{thm}
An element $a$ of commutative ring $A$ lies in the radical of an
ideal $\mathfrak{a}$ if and only if it lies in every prime ideal
containing $\mathfrak{a}$.
\end{thm}
Since we may assume that $f\notin\Delta(f)$, it suffices to show
the following.

\begin{prop}
If $f\in\mathbb{R}\left[\left[x_{1},x_{2},\ldots,x_{n}\right]\right]$
is a real power series with vanishing constant term and $\mathfrak{a}$
is a prime ideal of $\mathbb{R}\left[\left[x_{1},x_{2},\ldots,x_{n}\right]\right]$
containing the Jacobian ideal $\Delta(f)$, then $f\in\mathfrak{a}$.
\end{prop}
We shall approach this by first rephrasing in terms of derivations.

\subsubsection{Derivations}

\begin{note*}
All derivations here we presume to be \emph{continuous} derivations,
where continuity is with respect to the topology generated by powers
of the maximal ideal. Occasional mention will be made of this (usually
in the statements of lemmas, propositions, etc.), but the reader should
remember that the assumption of continuity is universal in what follows.
\end{note*}
The derivations $D:\mathbb{R}\left[\left[x_{1},x_{2},\ldots,x_{n}\right]\right]\rightarrow\mathbb{R}\left[\left[x_{1},x_{2},\ldots,x_{n}\right]\right]$
are of the form \[
\sum_{i=1}^{n}g_{i}(x)\frac{\partial}{\partial x_{i}}\]
 where $g_{i}(x)\in\mathbb{R}\left[\left[x_{1},x_{2},\ldots,x_{n}\right]\right]$,
$i=1,\ldots,n$. Indeed, if $D$ is a derivation, we let $g_{i}=D(x_{i})$.

Let $\mathfrak{a}\subset\mathbb{R}\left[\left[x_{1},x_{2},\ldots,x_{n}\right]\right]$
be an ideal, $A=\mathbb{R}\left[\left[x_{1},x_{2},\ldots,x_{n}\right]\right]/\mathfrak{a}$.
Then if $D:A\rightarrow A$ is a derivation, we can compose with the
quotient map to obtain a derivation $D^{\prime}:\mathbb{R}\left[\left[x_{1},x_{2},\ldots,x_{n}\right]\right]\rightarrow A$,
and by choosing a lift $g_{i}$ for each $D^{\prime}(x_{i})$ we obtain
a derivation on $\mathbb{R}\left[\left[x_{1},x_{2},\ldots,x_{n}\right]\right]$
which is a lift of the derivation on $A$. The following is then immediate.

\begin{lem}
Let $f$ be an element of the ring of formal power series $\mathbb{R}\left[\left[x_{1},x_{2},\ldots,x_{n}\right]\right]$
and let $\mathfrak{a}$ be an ideal of $\mathbb{R}\left[\left[x_{1},x_{2},\ldots,x_{n}\right]\right]$
containing the Jacobian ideal $\Delta(f)$ of $f$. Let $A=\mathbb{R}\left[\left[x_{1},x_{2},\ldots,x_{n}\right]\right]/\mathfrak{a}$
and $a=\left[f\right]\in A$. Then if $D:A\rightarrow A$ is a continuous
derivation, $D\left(a\right)=0$.
\end{lem}
We then obtain Proposition 1 from the following more general statement.

\begin{prop}
Let $\mathfrak{a}$ be a prime ideal in $\mathbb{R}\left[\left[x_{1},x_{2},\ldots,x_{n}\right]\right]$,
let $A$ be the quotient of $\mathbb{R}\left[\left[x_{1},x_{2},\ldots,x_{n}\right]\right]$
by $\mathfrak{a}$. If $a\in A$ satisfies\\
$D(a)=0$ for all continuous derivations $D:A\rightarrow A$, then
$a=0$.
\end{prop}
We shall make use of the fact that $\mathfrak{a}$ is a prime ideal
by passing to the quotient field of $A$.

\subsubsection{Quotient Field Results}

In this section we shall let $\mathfrak{a}$ be a prime ideal in $\mathbb{R}\left[\left[x_{1},x_{2},\ldots,x_{n}\right]\right]$,
$A$ be the quotient of $\mathbb{R}\left[\left[x_{1},x_{2},\ldots,x_{n}\right]\right]$
by $\mathfrak{a}$ and $a$ an element of $A$ with the property that
for all derivations $D:A\rightarrow A$, $D(a)=0$.

Since $\mathfrak{a}$ is prime, $A$ is an integral domain and we
let $F$ be its field of fractions.

\begin{lem}
Let $D:F\rightarrow F$ be a derivation of $F$. Then $D(a)=0$.
\end{lem}
\begin{proof}
Let $D:F\rightarrow F$ be a derivation of $F$. Applying $D$ to
$\left[x_{i}\right]$ we get elements $b_{i}$ of $F$. These $n$
elements are fractions $\frac{c_{i}}{d_{i}}$, $c_{i},d_{i}\in A$
and if we multiply $D$ by $\prod d_{i}$ we get a derivation $D^{\prime}$
with $D^{\prime}\left(\left[x_{i}\right]\right)\in A$. It follows
that $D^{\prime}:A\rightarrow A$, and so $D^{\prime}(a)=0$. Then
$D(a)=\frac{1}{\prod d_{i}}D^{\prime}(a)=0$.
\end{proof}
Now let $B=\mathbb{R}[a]$ be the smallest $\mathbb{R}$-algebra (with
identity) in $A$ containing $a$. We let $B[x_{1}]$ be the smallest
$\mathbb{R}$-algebra (with identity) in $A$ containing $B$ and
the class of $x_{1}$ in $A$. We let $B[[x_{1}]]$ be the closure
of $B[x_{1}]$ in $A$ with respect to the topology defined by the
powers of the maximal ideal of $A$. (This works well since $A$ is
a complete local ring.) We can then similarly form $B[[x_{1}]][x_{2}]$,
$B[[x_{1},x_{2}]]=B[[x_{1}]][[x_{2}]]$, etc. This gives us a chain
of inclusions \begin{eqnarray*}
B\subset B[x_{1}]\subset B[[x_{1}]]\subset B[[x_{1}]][x_{2}]\subset B[[x_{1},x_{2}]]\subset\\
\cdots\subset B[[x_{1},\ldots,x_{n}]] & = & A.\end{eqnarray*}
$K$ be the smallest field in $F$ containing $a$. We wish to show
that every derivation $D:K\rightarrow K$ extends to a derivation
$D':F\rightarrow F$. Note that the quotient fields of the rings in
the chain of inclusions above give a chain of inclusions of fields
\begin{eqnarray*}
K\subset K(x_{1})\subset K((x_{1}))\subset K((x_{1}))(x_{2})\subset K((x_{1},x_{2}))\subset\\
\cdots\subset K((x_{1},\ldots,x_{n})) & = & F.\end{eqnarray*}
 We shall extend $D$ step by step along this chain of fields. Each
step in the chain is of one of two forms. First, it is an extension
by a single element, \emph{e.g.}, from $K$ to $K(x_{1})$. Derivations
extend in these cases since such an extension is either algebraic
or transcendental. (See \cite{L}, p. 370.) Second, it is an extension
from the quotient field of some subring to the quotient field of the
closure of that subring. In this case, note first that a derivation
will be a multiple of a derivation on the subring, and it suffices
to extend the derivation on the subring to the closure of the subring.
However, since we are working in a complete local ring, we can obtain
the value of the extended derivation on a limit of elements of the
subring by taking the limit of the derivations of the elements of
the subring. This shows the following.

\begin{lem}
Every derivation $D:K\rightarrow K$ extends to a derivation $D':F\rightarrow F$.
\end{lem}
Thus

\begin{lem}
For every derivation $D:K\rightarrow K$, $D(a)=0$.
\end{lem}
We now apply the following (which is Proposition 5.2 in Chapter VIII,
\S 5, of \cite{L}).

\begin{prop}
A finitely generated extension $k(x)$ over $k$ is separable algebraic
if and only if every derivation $D$ of $k(x)$ which is trivial on
$k$ is trivial on $k(x)$.
\end{prop}
Applying this with $k=\mathbb{R}$ and $k(x)=K$, we see that $a$
is algebraic over $\mathbb{R}$. Thus $a$ satisfies some equation
\[
c_{n}a^{n}+c_{n-1}a^{n-1}+\cdots+c_{0}=0.\]
 However, $A$ still has a point; evaluation at the origin gives an
algebra map $A\to\R$. We can evaluate at the origin to obtain $c_{0}=0$.
Since $A$ is an integral domain, if $a\ne0$ we then have \[
c_{n}a^{n-1}+c_{n-1}a^{n-2}+\cdots+c_{1}=0,\]
 an impossibility if we started with the irreducible polynomial satisfied
by $a$. Hence $a$ must be zero, and we have proved Theorem \ref{thm:nilpotency}.

\subsection{The Family of Examples}

Let\[
A_{n}=\R\left[\left[x_{1},\ldots,x_{n}\right]\right]/\left(\frac{\partial f_{n}}{\partial x_{1}},\ldots,\frac{\partial f_{n}}{\partial x_{n}}\right),\]
where $f_{n}=\sum_{i=1}^{n}x_{i}^{3n-1}+\left(x_{1}\cdot\cdots\cdot x_{n}\right)^{3}$.
Let $a_{n}=\left[f_{n}\right]\in A_{n}$. Clearly, $a_{n}\in\left(A_{n}\right)_{\flat}$.
These provide a family of examples of the nilpotency we have just
proved. (These examples, and the following result, are mentioned at
the end of \cite{BS}.)

\begin{lem}
$\left(a_{n}\right)^{n-1}\ne0$ and $\left(a_{n}\right)^{n}=0$.
\end{lem}
\begin{proof}
The ideal is generated by the $n$ series (polynomials, actually)
$\left(3n-1\right)x_{i}^{3n-2}+\frac{3}{x_{i}}\left(x_{1}\cdot\cdots\cdot x_{n}\right)^{3}$,
$i=1,\ldots,n$. Thus \[
x_{i}^{3n-1}\equiv-\frac{3}{3n-1}\left(x_{1}\cdot\cdots\cdot x_{n}\right)^{3},\]
so \[
f_{n}\equiv\left(-\frac{3n}{3n-1}+1\right)\left(x_{1}\cdot\cdots\cdot x_{n}\right)^{3}=-\frac{1}{3n-1}\left(x_{1}\cdot\cdots\cdot x_{n}\right)^{3}.\]
Thus it suffices to show that both\[
\left(x_{1}\cdot\cdots\cdot x_{n}\right)^{3\left(n-1\right)}\not\equiv0\]
 and \[
\left(x_{1}\cdot\cdots\cdot x_{n}\right)^{3n}\equiv0.\]

To see the first, note that $\frac{\partial^{\left|\alpha\right|}}{\partial x_{1}^{\alpha_{1}}\cdots\partial x_{n}^{\alpha_{n}}}\left(g_{i}\right)\left(0\right)=0$
if all $\alpha_{j}\le3n-3$. By the product rule, the same is true
for anything in the ideal. However,\[
\frac{\partial^{n\left(3n-3\right)}}{\partial x_{1}^{3n-3}\cdots\partial x_{n}^{3n-3}}\left(\left(x_{1}\cdot\cdots\cdot x_{n}\right)^{3\left(n-1\right)}\right)\left(0\right)\ne0.\]

To see the second, note that\begin{eqnarray*}
\left(x_{1}\cdot\cdots\cdot x_{n}\right)^{3n-1} & = & x_{1}^{3n-1}\cdot\cdots\cdot x_{n}^{3n-1}\\
 & \equiv & \prod_{i=1}^{n}\left(-\frac{3}{3n-1}\right)\left(x_{1}\cdot\cdots\cdot x_{n}\right)^{3}\\
 & = & \left(-\frac{3}{3n-1}\right)^{n}\left(x_{1}\cdot\cdots\cdot x_{n}\right),^{3n}\end{eqnarray*}
so\[
0\equiv\left(x_{1}\cdot\cdots\cdot x_{n}\right)^{3n-1}\left(1-\left(-\frac{3}{3n-1}\right)^{n}\left(x_{1}\cdot\cdots\cdot x_{n}\right)^{3}\right).\]
The last factor on the right is an invertible power series, so\[
0\equiv\left(x_{1}\cdot\cdots\cdot x_{n}\right)^{3n-1},\]
and thus\[
0\equiv\left(x_{1}\cdot\cdots\cdot x_{n}\right)^{3n}.\]

\end{proof}
It is worthwhile noting here a way of representing elements of $A_{n}$.
Since $0\equiv\left(x_{1}\cdot\cdots\cdot x_{n}\right)^{3n}$, $x_{i}^{\left(3n-2\right)n}\equiv0$
for all $i$. Thus, by Taylor's theorem, any function $f\left(x_{1},\ldots,x_{n}\right)$
may be written\begin{eqnarray*}
f\left(x_{1},\ldots,x_{n}\right) & = & \tilde{P}\left(x_{1},\ldots,x_{n}\right)+\sum_{i=1}^{n}x_{i}^{\left(3n-2\right)n}R_{i}\left(x_{1},\ldots,x_{n}\right)\\
 & \equiv & \tilde{P}\left(x_{1},\ldots,x_{n}\right),\end{eqnarray*}
where $\tilde{P}\left(x_{1},\ldots,x_{n}\right)$ is a polynomial
of degree at most $\left(3n-2\right)n-1$ in each variable. Next,
any term in $\tilde{P}\left(x_{1},\ldots,x_{n}\right)$ having degree
greater than or equal to $3n-2$ in some $x_{i}$ can be replaced
(modulo the ideal) with one of degree $3n-4$ less in $x_{i}$ but
of total degree one greater. Repeating this process, each term either
becomes zero (having achieved too great a total degree) or becomes
of degree less than $3n-2$ in each $x_{i}$. Thus\[
f\left(x_{1},\ldots,x_{n}\right)\equiv P\left(x_{1},\ldots,x_{n}\right),\]
where $P\left(x_{1},\ldots,x_{n}\right)$ is a polynomial of degree
less than or equal to $3n-3$ in each of the variables $x_{i}$. Since
each of the generators of the ideal has a term of degree $3n-2$ in
some variable, the polynomial $P$ is uniquely determined by $\left[f\right]\in A_{n}$.

\subsection{The Second Splitting}

We now want to take these examples and, as we did before, split the
singularity into a finite number of points.

Let $a:\R\to\R$ be a smooth, monotone function vanishing to infinite
order at $0\in\R$. Let\[
h_{i}\left(x_{1},\ldots,x_{n},t\right)=\left(3n-1\right)\prod_{k=1}^{3n-2}\left(x_{i}-ka\left(t\right)\right)+\frac{3}{x_{i}}\left(x_{1}\cdot\cdots\cdot x_{n}\right)^{3}.\]
For fixed $t\ne0$, we will examine the common zeroes of the $h_{i}$.

Let $X_{i}=x_{i}/a$. Then\[
h_{i}=\left(3n-1\right)a^{3n-2}\prod_{k=1}^{3n-2}\left(X_{i}-k\right)+a^{3n-1}\frac{3}{X_{i}}\left(X_{1}\cdot\cdots\cdot X_{n}\right)^{3}.\]
Let\[
H_{i}=\frac{h_{i}}{a^{3n-2}}=\left(3n-1\right)\prod_{k=1}^{3n-2}\left(X_{i}-k\right)+a\frac{3}{X_{i}}\left(X_{1}\cdot\cdots\cdot X_{n}\right)^{3}.\]
At $t=0$, the common zeroes of the $H_{i}$ are $X_{i}=k_{i}$, $1\le k_{i}\le3n-2$.
Moreover, these are transverse intersections, so the common zeroes
will depend analytically on $a$. Since we have as many zeroes as
the dimension of $A_{n}$, this is the splitting we will use.

It is easy to generalize the representation of elements of $A_{n}$
described above to representation of functions $f\left(x_{1},\ldots,x_{n},t\right)\mod\left(h_{1},\ldots,h_{n}\right)$.
We can view such a function $f$ as an element of $A_{n}$ depending
on a parameter $t$. The reduction of $f$ to the polynomial $P$
depends smoothly on $t$, and we obtain\[
f=P\quad\mod\left(g_{1},\ldots,g_{n}\right),\]
where the coefficients of $P$ depend on $t$. That is, we write\[
f=P+\sum_{i=1}^{n}b_{i}g_{i}\]
for some elements $b_{i}$. We then have \[
f=P+\sum_{i=1}^{n}b_{i}h_{i}+F_{1}a.\]
Applying the same procedure to $F_{1}$, we obtain\[
f=P+\sum_{i=1}^{n}b_{i}h_{i}+F_{2}a^{2}\]
(for a different polynomial $P$). Proceeding in this fashion, we
obtain the following.

\begin{prop}
For all $N\in\mathbb{N}$ and for all $C^{\infty}$ functions $f=f\left(x_{1},\ldots,x_{n},t\right)$
defined in a sufficiently small neighborhood of the origin, there
is a polynomial $P\left(x_{1},\ldots,x_{n}\right)$ with degree at
most $3n-3$ in each of the variables $x_{i}$ and coefficients smooth
functions of $t$ such that\[
f\left(x_{1},\ldots,x_{n},t\right)\equiv P\left(x_{1},\ldots,x_{n}\right)\quad\mod\left(h_{1},\ldots,h_{n},a^{N}\right).\]

\end{prop}

\subsection{The Reasoning}

We now consider a representative $f$ of \[
\psi\left(\left[\left(x_{1}\cdot\cdots\cdot x_{n}\right)^{3}\right]\right)\]
 as an element of $C^{\infty}\left(U\times\R\right)$ for a small
neighborhood $U$ of the origin in $\R^{n}$, modulo the closed ideal
$I$ which is generated by $h_{1},\ldots,h_{n}$ and functions of
$t$ vanishing at $t=0$ and $t=\frac{1}{n},n\in\mathbb{N}$. $\left[\left(x_{1}\cdot\cdots\cdot x_{n}\right)^{3}\right]\in\left(C^{\infty}\left(U\times\R\right)/I\right)_{\flat}$,
so $f$ is well-defined modulo $I$. Moreover, $f$ will vanish at
all the points in $Z\left(I\right)$, by our assumptions on $\psi$.

Now there are unique polynomials\[
\sum_{i_{1},\ldots,i_{n}=0}^{3n-3}b_{i_{1}\cdots i_{n}}x_{1}^{i_{1}}\cdot\cdots\cdot x_{n}^{i_{n}}\]
interpolating values at the points\[
\left\{ \left(k_{1},\ldots,k_{n}\right):k_{i}\in\mathbb{Z},1\le k_{i}\le3n-2\right\} ,\]
the intersection points in the coordinates $X_{i}$ at $a=0$. (This
can be seen by first interpolating polynomials in $x_{1}$ for fixed
$\left(k_{2},\ldots,k_{n}\right)$, then interpolating polynomials
in $x_{2}$ (whose coefficients are polynomials in $x_{1}$) for fixed
$\left(k_{3},\ldots,k_{n}\right)$, etc.) In particular, the matrix
defining the linear map from the space of values to the space of such
polynomials in invertible. Since the zeroes are analytic in $a$,
for small $t$ the matrix is also invertible, and we have the same
interpolation. If we change coordinates back to the $x_{i}$, we have
the same sort of interpolation problem, but the determinant of the
matrix will vanish to some finite order $a^{K_{0}}$.

Now taking $N\gg K_{0}$, write\[
f=P+a^{N}f_{N}.\]
For $a\ne0$, the interpolating polynomial vanishes, so\[
0=P+a^{N-K_{0}}Q\]
 for some polynomial $Q$ depending smoothly on $t$. Now taking the
limit as $t\to0$, we find $\lim_{t\to0}P=0$, and hence at $t=0$,
$f\equiv0\quad\mod\left(g_{1},\ldots,g_{n}\right)$ as desired.

This completes the proof of the main theorem.

\end{document}